\documentclass[a4paper, 12pt]{article}
\usepackage{amsmath,amsthm,amsfonts,amssymb,stmaryrd,amscd}
\usepackage[top=2cm, bottom=2.5cm, left=2.2cm, right=2.2cm]{geometry}
\usepackage[french, english]{babel}

\input xy
\xyoption{all}

\usepackage{color}
\usepackage{makeidx}
\usepackage{hyperref}
\usepackage{pifont}

\makeindex 

\numberwithin{equation}{section}

{\theoremstyle{definition}
\newtheorem{definition}{Definition}

\newtheorem{remark}[definition]{Remark}

\newtheorem{remarks}[definition]{Remarks}

}
\newtheorem{proposition}[definition]{Proposition}
\newtheorem{proposition-definition}[definition]{Proposition-Definition}

\newtheorem{theorem}[definition]{Theorem}
\newtheorem{corollary}[definition]{Corollary}

\newcommand{\R}{\mathbb{R}}
\newcommand{\N}{\mathbb{N}}

\newcommand{\cF}{\mathcal{F}}

\newcommand{\cL}{\mathcal{L}}

\newcommand{\cW}{\mathcal{W}}

\newcommand{\cK}{\mathcal{K}}
\newcommand{\cX}{\mathcal{X}}

\newcommand{\cU}{\mathcal{U}}
\newcommand{\cH}{\mathcal{H}}
\newcommand{\cM}{\mathcal{M}}

\newcommand{\resp}{{\it resp.}\/ }

\newcommand{\diam}{{{\mathrm{diam}}}}
\newcommand{\ie}{{\it i.e.}\/ }
\newcommand{\eg}{{\it e.g.}\/ }
\newcommand{\cf}{{\it cf.}\/ }

\newcommand{\gA}{\mathfrak{A}}

\newcommand{\gag}{\mathfrak{g}}

\def\gpd{\,\lower1pt\hbox{$\longrightarrow$}\hskip-.24in\raise2pt
             \hbox{$\longrightarrow$}\,}

\reversemarginpar

\usepackage{color}

\everymath{\displaystyle}

\begin{document}

\begin{center}
{\Large \bf  Stability of Lie groupoid C*-algebras}
\footnote{The authors were partially supported by ANR-14-CE25-0012-01.} 

\bigskip
{\sc by Claire Debord and Georges Skandalis}

\end{center}

{\footnotesize
Laboratoire de Math\'ematiques, UMR 6620 - CNRS
\vskip-4pt Universit\'e Blaise Pascal, Campus des C\'ezeaux, BP {\bf 80026}
\vskip-4pt F-63171 Aubi\`ere cedex, France
\vskip-4pt claire.debord@math.univ-bpclermont.fr

\vskip 2pt Universit\'e Paris Diderot, Sorbonne Paris Cit\'e
\vskip-4pt  Sorbonne Universit\'es, UPMC Paris 06, CNRS, IMJ-PRG
\vskip-4pt  UFR de Math\'ematiques, {\sc CP} {\bf 7012} - B\^atiment Sophie Germain 
\vskip-4pt  5 rue Thomas Mann, 75205 Paris CEDEX 13, France
\vskip-4pt skandalis@math.univ-paris-diderot.fr
}
\bigskip

\centerline{\large \bf Abstract}
In this paper we generalize a theorem of M. Hilsum and G. Skandalis stating that the $C^*$- algebra of any foliation of non zero dimension is stable. Precisely, we show that the C*-algebra of a Lie groupoid is stable whenever the groupoid has no orbit of dimension zero. We also prove an analogous theorem for singular foliations for which the holonomy groupoid as defined by I. Androulidakis and G. Skandalis is not Lie in general.

\section{Introduction: Statement of the theorem and the steps of the proof}

The aim of this paper is to generalize Theorem 1 of \cite{HilsSkStab} stating that the $C^*$-algebra of any foliation (of non zero dimension !) is stable.

\begin{theorem}\label{casdesgroupoides}
 Let $G$ be a Lie groupoid with $\sigma$-compact  $G^{(0)}$. Assume that at every $x\in G^{(0)}$ the anchor $\natural_x:\gag_x\to T_x G^{(0)}$ is nonzero.
Then $C^*(G)$ is stable. 
\end{theorem}

In other words, $C^*(G)$ is stable whenever $G$ has no orbit of dimension $0$.

The converse is also true if $G$ is $s$-connected. Indeed, if $G$ is $s$-connected and the anchor at $x$ is the zero map, then the orbit of $x$ is reduced to $x$. Therefore $C^*(G)$ has a character: the trivial representation of the group $G_x^x$.

Since the reduced $C^*$-algebra $C^*_r(G)$ of $G$ is a quotient of $C^*(G)$, it follows that it is also stable when $G$ has no orbit of dimension $0$.

Here however, the converse may fail for the reduced $C^*$-algebra: the reduced $C^*$-algebra of the group $PSL_2(\R)$ is stable! 
  
Our proof is not very different from the one of \cite{HilsSkStab} and based on Kasparov's stabilisation theorem (\cite{KaspJOT}). Note that, unlike in   \cite{HilsSkStab}, we do not assume the space $G^{(0)}$ to be compact - but this is actually a rather minor point.

\medskip 
The proof is as follows.

\begin{enumerate}
\item \label{step1} Let $x\in G^{(0)}$. There is a section $Y$ of the algebroid $\gA(G)$ whose image under the anchor is a vector field $X$ satisfying $X(x)\ne 0$. Taking a local exponentiation of $X$ we obtain a relatively compact open neighborhood $W$ diffeomorphic to $U\times \R$ where $X$ is proportional to the vector field along the $\R$ lines $\{u\}\times \R$. 

This step  will be clarified in section \ref{Annexe}

\medskip 
We thus choose a locally finite cover $(W_n)$  by relatively compact open subsets and diffeomorphisms $f_n:U_n\times \R\to W_n$ such that $W'_n=f_n(U'_n\times \R)$ cover $G^{(0)}$ with $U'_n$ relatively compact in $U_n$.  Let $p_n:W_n\to U_n$ be the composition of $f_n^{-1}$ with the projection $U_n\times \R\to U_n$. 

\item \label{step2}One may then construct a locally finite family of open subsets $V_j$ of $G^{(0)}$ such that: \begin{itemize}
\item Every $V_j$ sits in a $W_{n(j)}$ and its intersection with each line $f_{n(j)}(u\times \R)$ is an (open) interval. More precisely,  ${f_{n(j)}(V_j})$ is of the form $\{(x,t)\in U_{n(j)}\times \R;\ \varphi_j^-(x)<t<\varphi_j^+(x)\}$ where $\varphi_j^-,\varphi_j^+:U_{n(j)}\to \R$ are smooth and $\varphi_j^+-\varphi_j^-$ is nonnegative with compact support. 

\item The $\overline{V_j}$ are pairwise disjoint and locally finite: every compact subset of $M$ meets only finitely many $\overline {V_j}$'s.

\item For every $n$, the $p_n(V_j\cap W_n)$ cover $U'_n$: we have $U'_n\subset \bigcup_{j;\ n(j)=n}p_n(V_j\cap W_n) $.
\end{itemize}

The details of the constructions of the $V_j$'s are given in section \ref{lesVj}.

\item \label{step3} Let then $q_j$ be the characteristic function of $V_j$. We prove that $q_j$ is a multiplier of $C^*(G)$. By local finiteness, the characteristic function  $q=\sum q_j$ of $V=\bigcup V_j$  is also a multiplier of $C^*(G)$. 

See section section \ref{VjMult}.

\item \label{step4} We show that $qC^*(G)$ is a full Hilbert submodule of $C^*(G)$ (see corollary \ref{fullsub} - section \ref{section3.2}). 

\item \label{step5} Considering a natural diffeomorphism $V_j\simeq p_n(V_j)\times {]}0,1{[}$, it follows that  the Hilbert $C^*(G)$-modules $q_jC^*(G)$ and $qC^*(G)$ are stable. 

\item \label{step6}Using Kasparov's stabilisation Theorem (\cite{KaspJOT}), it follows that $C^*(G)$ is stable. 

This follows from corollary \ref{corfinal} - see section \ref{stability}.
\end{enumerate}

In section \ref{foliation}, we prove an analogous theorem for singular foliations in the sense of \cite{AndrSk1}. We prove:

\begin{theorem}\label{casdesfeuilletages}
The $C^*$-algebra of a singular foliation (as defined in \cite{AndrSk1}) which has no leaves reduced to a single point is stable.
\end{theorem}

The main steps of the proof are the same as for theorem \ref{casdesgroupoides}.  Using vector fields along the foliation, we construct the same small open subsets $V_j$. Note that, in proving that the characteristic functions of these $V_j$ are multipliers of $C^*(M,\cF)$, we chose to take a somewhat different  path in order to shed a new light to it. This led us to construct groupoid homomorphisms between singular foliation groupoids (section \ref{lexion}). Of course, we could have used the same kind of proof as for the Lie groupoid case.

\section{Geometric constructions}

\subsection{Nonzero vector fields in the algebroid}\label{Annexe}

Let $M$ be a smooth (open) manifold, $x$ a point of $M$, and let $X\in \cX(M)$ be a smooth vector field with compact support on  $M$ such that $X(x)\ne 0$. Denote by $\Psi_X=(\Psi_X^t)_{t\in \R}$ the flow of $X$. 
One can find a codimension one submanifold $U$ of $M$ and a neighborhood $I$ of $0$ in $\R$ such that the restriction of $\Psi_X$ to $U\times I$ is a diffeomorphism onto an open (tubular) neighborhood $W$ of $U$ in $M$. In other words, $U$ is a codimension one submanifold  of $M$ which contains $x$ and which is transverse to the integral curves of $X$. 

\smallskip If $G$ is a Lie groupoid and $Y$ a section with compact support of its Lie algebroid such that the vector field $X:=\natural(Y)$ doesn't vanish on a point $x\in G^{(0)}$, let $Z_Y$ be the associated right invariant vector field on $G$ and $\Psi_{Z_Y}$ its flow. We have $r\circ \Psi_{Z_Y}^t=\Psi^t\circ r$.

 Applying the construction above, one finds a codimension one submanifold $U$ of $G^{(0)}$ and a neighborhood $I$ of $0$ in $\R$ such that the composition map  $$U\times I \overset{\Psi_{Z_Y}}{\longrightarrow} G \overset{r}{\longrightarrow} G^{(0)}$$ is precisely the restriction of the flow $\Psi_X$ of $X$ to $U\times I$ and thus a diffeomorphism onto an open neighborhood $W$ of $x$ in $G^{(0)}$. Note that $s\circ \Psi_{Z_Y}$ is the projection $U\times I \rightarrow U$.

Now the following maps are diffeomorphisms : 
$$\begin{array}{ccc} G^U\times I & \rightarrow & G^W \\ (\gamma,t) & \mapsto & \Psi_{Z_Y}(r(\gamma),t)\gamma \end{array} \ \mbox{ and } \ \begin{array}{ccc} G^U_U\times I \times I & \rightarrow & G^W_W \\ (\gamma,t,\lambda) & \mapsto & \Psi_{Z_Y}(r(\gamma),t)\gamma \Psi_{Z_Y}(s(\gamma),\lambda)^{-1}\end{array} \ .$$

\subsection{Construction of the family \boldmath{$V_j$}}\label{lesVj}

In this section, we explain the construction of the $V_j$'s. 

The construction above yields a locally finite cover $(W_n)$ by relatively compact open subsets and diffeomorphisms $f_n:U_n\times \R\to W_n$ such that $W'_n=f_n(U'_n\times \R)$ cover $G^{(0)}$ with $U'_n$ relatively compact in $U_n$. We will often identify $W_n$ and $U_n\times \R$ under $f_n$.

Let $p_n:W_n\to U_n$ be the composition of $f_n^{-1}$ with the projection $U_n\times \R\to U_n$. As $(W_n)$ is locally finite and $G^{(0)}$ is $\sigma$-compact, the set of indices is countable; we identify it with $\N$.

\begin{itemize}
\item First, choose a (Riemannian) metric $d$ on $G^{(0)}$. 
\item For every $n\in \N$, let $\varepsilon_n>0$ be small enough so that, for every $m\in \N$ such that $W_n\cap W_m\ne \emptyset$, and every $x\in U'_m$, $d(f_m(x,-1),f_m(x,1))>\varepsilon_n$. 
Thus if $W_n\cap W_m\ne \emptyset$, for every $x\in U'_m$ : $\diam(f_m(\{x\}\times \R)) > \varepsilon_n$.
\item Let $n\in \N$. Assume that we have already constructed a finite set $J_n$ with a map, $j\mapsto n(j)$ from $J_n$ to $\N$  with $n(j)<n$ and a family of open sets $(V_j)_{j\in J_n}$   with pairwise disjoint closures satisfying, 
 $$\forall  \ell<n ,\ \ U'_\ell\subset \bigcup_{j;\ n(j)=\ell}p_\ell(V_j), \ \ \diam(V_j)\le \varepsilon _{n(j)}.$$
 By the diameter assumption, for every $x\in U'_n$, the set $f_n(\{x\}\times \R)$ is not contained in a $\overline{V_j}$ with $j\in J_n$; therefore $f_n(\{x\}\times \R)$ is not contained in $\bigcup_{j\in J_n}\overline{V_j}$ (by connectedness of $\R$).

 By compactness, we may then construct a finite cover $(A_\ell)_{\ell\in J'_n}$ of $\overline{U'_n}$ by open subsets of $U_n$ and open intervals $I_\ell$ in $\R$ such that $\overline{f_n({A_\ell}\times {I_\ell})}\cap \bigcup_{j\in J_n}\overline{V_j}=\emptyset$ and $\diam(f_n({A_\ell}\times I_\ell))\le\varepsilon _n$.
 
Replacing the $I_\ell$'s by smaller intervals, we now further assume that the $\overline{I_\ell}$ with $\ell\in J'_n$ are pairwise disjoint. One then constructs for each $\ell\in J'_n$ two smooth functions $\varphi_\ell^{\pm}:U_n\to I_\ell$ such that $\varphi_\ell^{+}\ge \varphi_\ell^{-}$ and such that $A_\ell=\{x\in U_n; \ \varphi_\ell^{-}(x)< \varphi_\ell^{+}(x)\}$. 
 
Let then $J_{n+1}$ be the disjoint union of $J_n$ end $J'_n$. For $\ell\in J'_n$, put $V_\ell=f_n(\{(x,t)\in U_n\times \R;\ \ \varphi_\ell^{-}(x)< t<\varphi_\ell^{+}(x)\})$   and define $n(\ell)=n$.

\item We will also use the diffeomorphism $h_j:A_j\times {]}0,1{[}\to V_j$ given by \label{carreau}
\begin{equation}\tag{\ding{71}}
 h_j(x,t)=f_{n(j)}\Big(x,(1-t)\varphi_j^{-}(x)+ t\varphi_\ell^{+}(x)\Big).
\end{equation}

\item We thus construct the family $(V_j)$ inductively. Every compact set meets only finitely many $W_n$'s  and therefore finitely many $\overline{V_j}$'s since $\overline{V_j}\subset W_{n(j)}$ and $\{j;\ n(j)=n\}$ is finite.
\end{itemize}

 \section{Hilbert modules}
 
 \subsection{Stable Hilbert modules; full Hilbert modules}
 
We start by briefly recalling some now classical general facts on Hilbert modules. The basic reference for them is \cite{KaspJOT}. See also e.g. \cite{WE, Lance}.
 
Let $A$ be a separable $C^*$-algebra. 

\begin{definition}
A separable Hilbert $A$-module $E$ is said to be\begin{description}
\item[Stable] if it is isomorphic to $\ell^2(\N)\otimes E$.
\item[Full] if the vector span of the set of products $\langle x|y \rangle$ with $x,y\in E$ is dense in $A$.
\end{description}
\end{definition}

We have:

\begin{proposition}\label{utation}
\begin{enumerate}
\item Let $J$ be a countable set and $(E_j)_{j\in J}$ a family of (separable) stable Hilbert $A$ modules, then $\bigoplus_{j\in J} E_j$ is stable.
\item A separable Hilbert $A$-module $E$ is stable if and only if the $C^*$-algebra $\cK(E)$ is stable (see \eg \cite[Facts 4.7]{DS1}). \hfill$\square$
\end{enumerate}
\end{proposition}

Let us also recall Kasparov's stabilization Theorem. Put $\cH_A=\ell^2(\N)\otimes A$.

\begin{theorem}[Kasparov's stabilization Theorem]
For every separable Hilbert $A$-module $E$, the Hilbert $A$-modules $\cH_A$ and $E\oplus \cH_A$ are isomorphic.
 \hfill$\square$
\end{theorem}

The following statement is an immediate consequence of Kasparov's stabilization theorem. We outline a proof for completeness.

\begin{corollary}\label{corfinal}
\begin{enumerate}
\item A separable full and stable Hilbert $A$-module is isomorphic to $\cH_A$ (\cf \cite[Prop. 7.4, p. 73]{Lance}).
\item If $E$ is a separable Hilbert $A$-module and $p\in \cL(E)$ is a projection such that $pE$ is full and stable, then $E$  is isomorphic to $\cH_A$.
\end{enumerate}
\begin{proof}
\begin{enumerate}
\item  Put $B=\cK(E)$ and $E^*=\cK(E,A)$ considered as a Hilbert $B$ module. Since $E$ is full, we find $\cK(E^*)=E^*\otimes _BE=A$. Since $E$ is stable, the Hilbert $B$-module $B$ is isomorphic to $\cH_B$, therefore, using Kasparov's stabilization theorem, we find an isomorphism $u:B\to B\oplus (\ell^2(\N)\otimes E^*)$. We thus obtain an isomorphism of Hilbert $A$-modules $u\otimes_B1_{E}$ between $E=B\otimes_BE$ and $E\oplus (\ell^2(\N)\otimes E^*\otimes_B E)\simeq \cH_A$.

\item Write $E=pE\oplus (1-p)E\simeq \cH_A\oplus (1-p)E\simeq \cH_A$.
\qedhere
\end{enumerate}
\end{proof}
\end{corollary}

\subsection{Hilbert module associated with a transverse map}\label{section3.2}

We now recall the rather well known construction of generalized smooth transverse maps $M\to G$ from a manifold $M$ to a smooth groupoid $G$. This construction is useful in steps \ref{step4} and \ref{step5} of the proof of our theorem as explained in section 1. Note that we actually only use it for a (locally closed) submanifold $M\subset G^{(0)}$  transverse to $G$. 

\begin{definition}
Let $M$ be a smooth manifold and $G$ a smooth groupoid with algebroid $\gA$ and anchor $\natural$. A smooth map $f:M\to G^{(0)}$ is said to be \emph{transverse} to $G$ if for every $x\in M$, $df_x(T_xM)+\natural _{f(x)}\gA_{f(x)}=T_{f(x)}G^{(0)}$.
\end{definition}

The \emph{graph} of a smooth map $f:M\to G^{(0)}$ is the set $\Gamma^f=M\times _{G^{(0)}}G=\{(x,\gamma)\in M\times G;\ f(x)=r(\gamma)\}$.

Equivalently, $f$ is transverse if and only if the \emph{source map} $s_f$ is a submersion, where $s_f:\Gamma^f \to  G^{(0)}$ is the map $(x,\gamma) \mapsto  s(\gamma)$.

We will in fact only use the case where $f:M\to G^{(0)}$ is the inclusion of a (locally closed) submanifold $M$ of $G^{(0)}$. Then, $\Gamma^f=G^M=\{\gamma\in G;\ r(\gamma)\in M\}$.

\medskip For the sake of completeness, we recall also the notion of a \emph{generalized transverse map.}

Let $G$ be a smooth groupoid, $M$ a smooth manifold. A generalized morphism  $f:M\to G$ is given either by:
 \begin{description}
\item[A cocycle]: An open cover $(U_i)_{i\in I}$ of $M$, smooth maps $f_i:U_i\to G^{(0)}$ and smooth maps $f_{ij}:U_i\cap U_j\to G$ satisfying $r\circ f_{ij}=f_i$, $s\circ f_{ij}=f_j$ and, for all $x\in U_i\cap U_j\cap U_k$, $f_{ik}(x)=f_{ij}(x)f_{jk}(x)$.

\item[The graph of $f$]: A set $\Gamma^f$ which is a $G$-principal bundle over $M$. We therefore are given maps $r_f:\Gamma^f\to M$ which is a smooth surjective submersion and $s_f:\Gamma^f\to G^{(0)}$  with a right action of $G$ which is a smooth map $(x,\gamma)\mapsto x\gamma$ from $\Gamma^f\times _{G^{(0)}}G=\{(x,\gamma)\in \Gamma^f\times G;\ s_f(x)=r(\gamma)\}$ to $\Gamma^f$. This action is assumed to be proper and free with quotient $M$:
\begin{description}
\item[Proper] means that the map $(x,\gamma)\mapsto (x,x\gamma)$ is proper from $\Gamma^f\times _{G^{(0)}}G$ to $ \Gamma^f\times \Gamma^f$.
\item [Free] means that the map $(x,\gamma)\mapsto (x,x\gamma)$ is injective.
\item[The quotient is $M$] means that for a pair $(x,y)\in \Gamma^f\times \Gamma^f$ there exists $\gamma\in G$ with $x\gamma =y$ if and only $r_f(x)=r_f(y)$. Note that by freeness, this $\gamma$ is unique. We will denote it by $x^{-1}y$.
\end{description}
These three conditions altogether mean that the map $(x,\gamma)\mapsto (x,x\gamma)$ is a diffeomorphism from $\Gamma^f\times _{G^{(0)}}G$ onto $\Gamma^f\times _M\Gamma^f=\{(x,y)\in \Gamma^f\times \Gamma^f;\ r_f(x)=r_f(y)\}$ - whose inverse is $(x,y)\mapsto (x,x^{-1}y)$.
\end{description}
Given a cocycle $(f_i)$, $(f_{i,j})$, we obtain the graph $\Gamma^f$ by gluing the graphs  $\Gamma^{f_i}$ thanks to the $f_{i,j}$'s. Conversely, we obtain a cocycle out of $\Gamma^f$ by means of local sections of the submersion $r_f:\Gamma^f\to M$.

The generalized morphism $f$ is \emph{transverse} to $G$ (or a submersion from $M$ to the `bad manifold' $G^{(0)}/G$) if $s_f:\Gamma^f\to G^{(0)}$ is a submersion. This is equivalent to saying in the cocycle vision that the maps $f_{i}$ are transverse to $G$.

If $f$ is transverse to $G$,  the map $(x,\gamma)\mapsto \gamma$ is a submersion from $\Gamma^f\times _{G^{(0)}}G$ into $G$, whence the map $(x,y)\mapsto x^{-1}y$ is a submersion from $\Gamma^f\times _M\Gamma^f$ to $G$

It then defines a Hilbert-$C^*(G)$-module $C^*(\Gamma^f)$: this is the completion of $C_c(\Gamma^f)$ with respect to the $C^*(G)$-valued inner product given (using $1/2$-densities) by $$\langle g|h\rangle(\gamma)=\int_{x\in \Gamma^f;\ s_f(x)=r(\gamma)} \overline g(x)h(x\gamma).$$

Let us state the following easy fact which is important for our constructions:

\begin{proposition}
 Let $f:M\to G$ be a generalized \emph{transverse} morphism with graph $\Gamma^f$. The module $C^*(\Gamma^f)$ is full if and only if $s_f:\Gamma^f\to G^{(0)}$ is onto. 
\begin{proof}
The image of the submersion $s_f:\Gamma^f\to G^{(0)}$ is an open subset $U$ of $G^{(0)}$. Using the action of $G$ on $\Gamma^f$, we deduce that $U$ is saturated in $G^{(0)}$ - i.e. if $r(\gamma)\in U$, there exists $x\in \Gamma^f$ with $r(\gamma)=s_f(x)$; then $s(\gamma)=s_f(x\gamma)\in U$.
 
The image of the map $(x,y)\mapsto x^{-1}y$ from (defined on  $\Gamma^f\times _M\Gamma^f$) is $G_U^U$.
 
It follows that the inner products $\langle g|h\rangle$ with $g,h\in C_c^\infty(\Gamma^f)$, span a dense subset of $C_c^\infty(G_U^U)$.
 
 If $U=G^{(0)}$, these scalar products span a dense subspace in $C^*(G)$; if $U\ne G^{(0)}$, they all sit in the kernel of a regular representation associated to any point $x\not\in U$.
\end{proof}
\end{proposition}

Note that if $M$ is a submanifold of $G^{(0)}$, this condition means that $M$ meets all the $G$ orbits. 

Let us state this proposition in the precise way we will need to use it:

\begin{corollary}\label{fullsub}
Let $G$ be smooth and $V$ an open subset of $G^{(0)}$. Consider $C_0(V)$ as sitting in $C_0(G^{(0)})$ and therefore in the multiplier algebra of $C^*(G)$. Assume that every orbit of $G$ has a nonempty intersection with $V$. Then $C_0(V)C^*(G)$ is a full submodule of $C^*(G)$.\hfill$\square$
\end{corollary}

\subsection{The characteristic function of \boldmath{$V_j$} is a projection}\label{VjMult}

We prove here that the characteristic function $q_j$ of $V_j$ is a multiplier of $C^*(G)$, \ie that $C_0(V_j)C^*(G)$ is \emph{orthocomplemented} in $C^*(G)$.

Indeed, put $n=n(j)$. We may write $q_j=\vartheta_j q_j\vartheta_j$ where $\vartheta_j$ is a smooth real valued function on $G^{(0)}$ with support in $W_n$  which is equal to $1$ on $\overline V_j$. It is then enough to prove that $q_j\in \cL(C^*(G^{W_{n}}))$ viewing  $\vartheta_j$ as an element $T_{\vartheta_j}\in \cL(C^*(G);C^*(G^{W_{n}}))$ and write $q_j=T_{\vartheta_j}^*q_jT_{\vartheta_j}$.

Using the  identification between $W_n$ and $U_n\times \R$ coming from the diffeomorphism  $f_n$, we may write $C^*(G^{W_{n}})  =(C_0(U_n)\otimes L^2(\R))\otimes_{C_0(U_n)}C^*(G^{U_{n}})$. 

The characteristic function of $V_j$ is given by a ($*$-)strongly continuous map from $U_{n}$ to $\cL(L^2(\R))$, and therefore is an element $Q_j \in \cL(C_0(U_n)\otimes L^2(\R))$. Therefore $q_j=Q_j\otimes_{C_0(U_n)}1$ is in $\cL(C^*(G^{W_{n}}))$. 

Note also that  $Q_j$ is the $*$-strong limit of a sequence $(\theta_k)$ of smooth functions on $U_n\times \R$: put $\theta_k(x,t)=\phi(k(t-\varphi^- (x)))\phi(k(\varphi^+ (x)-t))$ where $\phi:\R\to [0,1]$ is a continuous function vanishing for $t\le 0$ and equal to $1$ for $t\ge 1$. It follows that the range of $q_j$ is (the closure of) $C_0(V_j)C^*(G)$. It is the Hilbert $C^*(G)$-module $C^*(G^{V_j})$ corresponding to the inclusion $V_j\to G^{(0)}$.

Now the projections $q_j$ are pairwise orthogonal and, by local finiteness, the sum $\sum q_j$ is strictly convergent. Indeed, the sum $\sum q_j\xi$ has only finitely many nonzero terms for $\xi\in C_c(G)$.

\subsection{Stability} \label{stability}

\begin{description}
\item[Stability of \boldmath{$q_jC^*(G)$}.] Using the diffeomorphism $h_j:A_j\times {]}0,1{[}\to V_j$ (see section \ref{carreau}, formula (\ding{71})), we deduce a diffeomorphism $G^{V_j}\simeq G^{A_j}\times {]}0,1{[}$, whence an isomorphism of Hilbert $C^*(G)$-modules $C^*(G^{V_j})\simeq C^*(G^{A_j})\otimes L^2({]}0,1{[})$. Therefore, the Hilbert $C^*$-module $q_jC^*(G)$ is stable.

\item[Stability of \boldmath{$qC^*(G)$}.] Since $q=\sum q_j$ ($*$-strong convergence), it follows that $qC^*(G)=\bigoplus_j q_jC^*(G)$. Therefore $qC^*(G)$ is stable too.

\item[Conclusion.] It then follows from corollary \ref{corfinal} that the Hilbert $C^*(G)$-module $C^*(G)$ is isomorphic to $\cH_{C^*(G)}$. The $C^*$-algebra $C^*(G)=\cK(C^*(G))=\cK(\cH_{C^*(G)})$ is stable. 
\end{description}
This ends the proof of theorem \ref{casdesgroupoides}.

\section{Stability of the C*-algebra of a singular foliation}\label{foliation}

\subsection{The holonomy groupoid of a singular foliation}

The C*-algebra of a singular foliation was defined in \cite{AndrSk1}. Let us briefly recall a few facts and constructions from \cite{AndrSk1}.

Recall that a foliation on a manifold $M$ is defined in \cite{AndrSk1} to be a (locally) finitely generated submodule $\cF$, stable by Lie brackets, of the $C^\infty(M)$-module $\cX_c(M)$ of smooth vector fields on $M$ with compact support.

A \emph{bi-submersion} of $\cF$ is the data of $(N,r_N,s_N)$ where $N$ is a smooth manifold, $r_N,s_N : N\rightarrow M$ are smooth submersions such that : 
$$r_N^{-1}(\cF)=s_N^{-1}(\cF) \mbox{ and } s_N^{-1}(\cF)=C_c^{\infty}(N;\mbox{ker } ds_N) + C_c^{\infty}(N;\mbox{ker } dr_N) \  \footnote{If $h:N \rightarrow M$ is a smooth submersion $h^{-1}(\cF)$ is the vector space generated by tangent vector fields $fZ$ where $f\in C_c^{\infty}(N)$ and $Z$ is a smooth tangent vector field on $N$ which is projectable by $dh$ and such that $dh(Z)$ belongs to $\cF$.}\!^) .$$
The \emph{inverse} of $(N,r_N,s_N)$ is $(N,s_N,r_N)$ and if $(T,r_T,s_T)$ is another bi-submersion for $\cF$ the \emph{composition}
is given by $(N,r_N,s_N)\circ (T,r_T,s_T):=(N \times_{s_N,r_T} T,r_N\circ p_N, s_T\circ p_T) \ ,$
where $p_N$ and $p_T$ are the natural projections respectively of $N \times_{s_N,r_T} T$ on $N$ and on $T$. \\
A \emph{morphism} from $(N,r_N,s_N)$ to $(T,r_T,s_T)$ is a smooth map $h:N\rightarrow T$ such that $s_T\circ h=s_N$ and $r_T\circ h=r_N$ and it is \emph{local} when it is defined only on an open subset of $N$. \\ 
Finally a bi-submersion can be \emph{restricted} : if $U$ is an open subset of $N$, $(U,r_U,s_U)$ is again a bi-submersion, where  $r_U$ and $s_U$ are the restriction of $r_N$ and $s_N$ to $U$.

\smallskip \noindent  For $x$ in $M$, we define the  \emph{fiber} of $\cF$ at $x$ to be the quotient $\cF_x=\cF/ I_x\cF$.
Let $\cX=(X_i)_{i\in \llbracket  1,n\rrbracket} \in \cF^n$ be such that $\cX_x=([X_i]_x)_{i\in \llbracket  1,n\rrbracket}$ is a basis of $\cF_x$. For any $\xi=(\xi_i)_{i\in \llbracket  1,n\rrbracket} \in \R^n$, we consider the vector field $X_{\xi}:= \sum_{i=1}^n \xi_i X_i$ and we denote by $\Psi_{\xi}^s$ its flow at time $s$. We consider the two smooth submersions from $M \times \R^n$ to $M$ : 
$$(s_{\cX},r_{\cX}) :  M\times \R^n  \longrightarrow  M\times M \ ; \  (x,\xi)  \mapsto  (x,\Psi_{\xi}^1(x)) \ .$$

According to Proposition 2.10 and 3.11 of \cite{AndrSk1}, one can find an open neighborhood $W$ of $(x,0)$ in $M\times \R^n$ such that $(W,r_W,s_W)$ is a bi-submersion, where the map $r_W$ and $s_W$ are the restriction to $W$ of the maps $r_{\cX}$ and $s_{\cX}$ defined above. Such a bi-submersion is called a \emph{path holonomy bi-submersion minimal at} $x$.

\smallskip \noindent Notice, that any restriction around $(x,0)$ of a path holonomy bi-submersion minimal at $x$, is again a path holonomy bi-submersion minimal at $x$.
  
\smallskip \noindent Let $\cU=(U_i,r_i,s_i)_{i\in I}$ be a family of bi-submersions of $\cF$. A bi-submersion  $(U,r_U,s_U)$ is \emph{adapted} to $\cU$ if for any $u\in U$ there is a local morphism around $u$ from $U$ to a $U_i$ for some $i$. 

\smallskip \noindent A family $\cU=(U_i,r_i,s_i)_{i\in I}$ of bi-submersions of $\cF$ which satisfies $M=\underset{i\in I}{\bigcup} s_i(U_i)$ and the inverse of any element in $\cU$ is adapted to $\cU$ together with the composition of any two elements of $\cU$ is an \emph{atlas}.

The \emph{path holonomy atlas} is the family of bi-submersions of $\cF$ generated by the path holonomy bi-submersions.

\smallskip \noindent The \emph{groupoid of the atlas} $\cU$ is the quotient $G(\cU)=\sqcup_{i\in I} U_i/\sim$ where $U_i\ni u \sim v\in U_j$ if and only if there is a local morphism from $U_i$ to $U_j$ sending $u$ on $v$. When $(U,r_U,s_U)$ belongs to $\cU$ and $u\in U$, let us denote by $[U,r_U,s_U]_u$  its image in $G(\cU)$. The structural morphisms of $G(\cU)$ are given by :
\begin{description}
\item[source and range :] ${\bf s}([U,r_U,s_U]_u)=s_U(u),   {\bf r}([U,r_U,s_U]_u)=r_U(u)$,
\item[inverse :] $[U,r_U,s_U]_u^{-1}=[U,s_U,r_U]_u$,
\item[product :]$[U,r_U,s_U]_u \cdot [V,r_V,s_V]_v=[U\times_{s_U,r_V}V,r_U\circ p_U, s_V\circ p_V]_{(u,v)} \mbox{ when } s_U(u)=r_V(v)$.
\end{description}

\noindent The groupoid of an atlas is endowed with the quotient topology which is quite bad, in particular the dimension of the fibers may change.

The \emph{holonomy groupoid} of $\cF$ is the groupoid of the path holonomy atlas.

\subsection{Subfoliations}\label{lexion}

A \emph{subfoliation} $\cF_1$ of a foliation $\cF_2$ is a submodule of $\cF_2$ which is a foliation \ie it is locally finitely generated and stable by Lie brackets.

In this section, we fix a foliation $\cF_2$ and a subfoliation $\cF_1$ of $\cF_2$.

\subsubsection{The atlas of compatible bi-submersions}

\begin{definition}
A bi-submersion $(U,r,s)$ of $\cF_1$ is said to be \emph{compatible} with $\cF_2$ if $r^{-1}(\cF_2)=s^{-1}(\cF_2)$.
\end{definition}

\begin{proposition}\label{ge}
\begin{enumerate}
\item For every $x_0\in M$, there is a  bi-submersion $(\cW,r,s)$ of $\cF_1$  compatible with $\cF_2$ such that $x_0\in s(\cW)$.
\item If $f:(U,r_U,s_U)\to (V,r_V,s_V)$ is a morphism of bi-submersions for $\cF_1$ and $(V,r_V,s_V)$ is compatible with $\cF_2$, then $(U,r_V,s_V)$ is compatible with $\cF_1$.
\item The bi-submersions of $\cF_1$ compatible with $\cF_2$ form an atlas.\footnote{To make it a set, take only bi-submersions defined on open subsets of $\R^N$ for all $N$. Note that a bi-submersion $(U,r,s)$ of $\cF_1$ is then adapted to this atlas if and only if $r^{-1}\cF_2=s^{-1}\cF_2$.}
\end{enumerate}
\begin{proof}
\begin{enumerate}
\item Fix  vector fields $X_1,\ldots ,X_n$ which generate $\cF_1$ in a neighborhood of $x_0$. For $y=(y_1,\ldots,y_n)\in \R^n$, put $\varphi_y=\exp (\sum y_iX_i)\in \exp \cF$. On $\R^n\times M$, put $s_0(y,x)=x$ and $t_0(y,x)=\varphi_y(x)$.  It is proved in  \cite[section 2.3]{AndrSk1}  that there is a neighborhood $\cW $ of $(0,x_0)$ in $\R^n\times M$ such that $(\cW ,t,s)$ is a bi-submersion for $\cF_1$ where $s$ and $t$ are the restrictions of $s_0$ and $t_0$.

Let $Y$ be the vector field on $\R^n\times M$ given by $Y(y,x)=(0,\sum y_iX_i)\in \R^n\times T_xM=T_{(y,x)}(\R^n\times M)$. Since $Y\in s^{-1}\cF_1\subset s^{-1}\cF_2$, it follows that $s^{-1}\cF_2$ is invariant under $\exp Y$. We find that $s^{-1}(\cF_2)=(s\circ \exp Y)^{-1}(\cF_2)$ as desired.

\item As $r_U=r_V\circ f$ and $s_U=s_V\circ f$, we find $r_U^{-1}(\cF)=f^{-1}(r_V^{-1}(\cF))=f^{-1}(s_V^{-1}(\cF))=s_U^{-1}(\cF)$.

\item Let $(U,r_U,s_U)$ and $(V,r_V,s_V)$ be bi-submersions (for $\cF_1$) compatible with $\cF_2$.\begin{enumerate}
\item Obviously, the inverse $(U,s_U,r_U)$ of $(U,r_U,s_U)$ is compatible with $\cF_2$.
\item Recall that the composition $(W,r_W,s_W)$ of $(U,r_U,s_U)$ and $(V,r_V,s_V)$ is constructed as follows: define $W=\{(u,v)\in U\times V;\ s_U(u)=r_V(v)\}$ and for $(u,v)\in W$,  set $s_W(u,v)=s_V(v)$ and $r_W(u,v)=r_U(u)$.  As $p_U:(u,v)\mapsto u$ and $p_V(u,v)\mapsto v$ are submersions on $W$, and $s_U\circ p_U=r_V\circ p_V$, we find 
\begin{eqnarray*}
r_W^{-1}(\cF_2)&=&p_U^{-1}(r_U^{-1}(\cF_2))=p_U^{-1}(s_U^{-1}(\cF_2))\\
&=&p_V^{-1}(r_V^{-1}(\cF_2))=p_V^{-1}(s_V^{-1}(\cF_2))\\
&=&s_W^{-1}(\cF_2).
\end{eqnarray*}\qedhere
\end{enumerate}

\end{enumerate}
\end{proof}
\end{proposition}

\begin{proposition}\label{us}
The composition of a bi-submersion of $\cF_1$ compatible with $\cF_2$ with a bi-submersion of $\cF_2$ is a bi-submersion of $\cF_2$.

\begin{proof} 
Let $(U,r_U,s_U)$ be a bi-submersion of $\cF_1$ compatible with $\cF_2$ and $(V,r_V,s_V)$ be a bi-submersion for $\cF_2$. As previously let $W=\{(u,v)\in U\times V;\ s_U(u)=r_V(v)\}$ and denote by $p_U$ and $p_V$ the projections on $U$ and $V$. We set $s_W=s_V\circ p_V$, $r_W=r_U\circ p_U$ and $\alpha=s_U\circ p_U=r_V\circ p_V$.

Since  $(U,r_U,s_U)$ is compatible with $\cF_2$ we have : $$r_W^{-1}(\cF_2)=\alpha^{-1}(\cF_2)=s_W^{-1}(\cF_2) $$
In particular $C_c^{\infty}(W;\ker dr_W)+C_c^{\infty}(W;\ker ds_W)\subset s_W^{-1}(\cF_2)$.

On the other hand $$\begin{array}{ccl} 
C_c^{\infty}(W;\ker d\alpha) & = & C_c^{\infty}(W;\ker dp_U) \oplus C_c^{\infty}(W;\ker dp_V) \\ & \\  s_W^{-1}(\cF_2) &= & p_V^{-1}(C_c^{\infty}(V;\ker ds_V) + C_c^{\infty}(V;\ker dr_V) )\\ & = & C_c^{\infty}(W;\ker ds_W) +C_c^{\infty}(W;\ker d\alpha) \\ & = &  C_c^{\infty}(W;\ker dp_U) + C_c^{\infty}(W;\ker dp_V) + C_c^{\infty}(W;\ker ds_W) \\ & \subset & C_c^{\infty}(W;\ker dr_W)+C_c^{\infty}(W;\ker ds_W)+C_c^{\infty}(W;\ker ds_W)
\end{array}$$
\end{proof}
\end{proposition}

\begin{remark}\label{rigerer}
Given bi-submersions $(U,r_U,s_U)$ and $(U',r_U',s_U')$ of $\cF_1$ and $(V,r_V,s_V)$ and $(V',r_V',s_V')$ of $\cF_2$ and morphisms $f:U\to U'$ and $g:V\to V'$, we obtain a morphism $(u,v)\mapsto (f(u),g(v))$ of bi-submersions $U\circ V\to U'\circ V'$.
\end{remark}

\subsubsection{Morphism on holonomy groupoids}

Let $G_i$ be the (path) holonomy groupoids of $\cF_i$ (\cite[Example 3.4.3 and Def. 3.5]{AndrSk1}). Let $\check G_1$ be the groupoid associated with the atlas of bi-submersions of $\cF_1$ compatible with $\cF_2$.  Let $\widehat G_2$ be the maximal holonomy groupoid of $\cF_2$ associated to the ``maximal atlas'' of all possible bi-submersions (\cite[Ex. 3.4.1]{AndrSk1}). 

Let $(U,r_U,s_U)$ in $\cU(\cF_1)$, $u\in U$ and $x=s_U(u)$. Choose a bi-submersion  $(V,r_V,s_V)$ and $v\in V$ representing the unit $x\in G_2^{(0)}$ that is :
\begin{itemize} 
\item $v\in V$ satisfies $s_V(v)=r_V(v)=x$ and 
\item There is an   identity bisection through $v$. This means that there is an open subset $W\subset M$ and a smooth map $j:W\to V$ which is a section for both $s$ and $r$.
\end{itemize}
According to proposition \ref{us}, the composition of $(U,r_U,s_U)$ and $(V,r_V,s_V)$ is a bi-submersion for $\cF_2$. Let then $U'$ be the open subset $s_U^{-1}(W)$ of $U$ and define $\psi :U'\to U\circ V$ by $\psi (u)=(u,j\circ s(u))$.

According to remark \ref{rigerer}, the class of $\psi(u)$ in $\widehat G_2$ only depends on the class of $u$ in $\check G_1$. Note also, that taking a bi-submersion  $(V'',r''_V,s''_V)$ and $v''\in V''$ representing the unit $r_U(u)\in G_2^{(0)}$, we find a bi-submersion $V''\circ U$ for $G_2$. As the class of $v$ (\resp $v''$) is the identity at $x$ (\resp $r_U(u)$), the classes of $(v'',u)\in V''\circ U$ (\resp $(u,v)\in U\circ V$) and $(v'',u,v)\in V''\circ U\circ V$ coincide.

Finally, for any bi-submersion $V_1$ for $\cF_1$, bi-submersion $V_2$ for $\cF_2$  the class $[(v_1,v_2)]$ of $(v_1,v_2)\in V_1\circ V_2$ is $\varphi([v_1])[v_2]$.

We thus constructed a ``smooth'' morphism of groupoids $\varphi:\check G_1\to \widehat G_2$. Smoothness means that for every $\gamma_1\in G_1$ there is a bi-submersion $U$ associated with $\check G_1$ such that $\gamma_1$ is the class of some element in $U$  and a smooth map $\psi$ from $U$ is a bi-submersion  $V$ for $\cF_2$ such that the diagram 
$$\xymatrix{U\ar[r]^{\psi}\ar[d]&V\ar[d]\\ \check G_1\ar[r]^\varphi& \widehat G_2}$$
commutes, where the vertical maps associate to elements of the bi-submersions their class in the corresponding groupoids.

According to \cite{DebCRAS}  and remark 3.13 of \cite{AndrSk1}, the $s$-fibers $\check G_{1,x}$ and $\widehat G_{2,x}$ are smooth manifolds. The restriction of $\varphi$ defines in particular a smooth map $\check G_{1,x}\to \widehat G_{2,x}$.

By definition of the path holonomy atlas, the corresponding groupoid $G_1$ is $s$-connected, since it is generated by any small neighborhood of the identity. More precisely, and thanks to \cite{DebCRAS}, the $s$-fibers  $G_{1,x}$ are connected manifolds. It follows that $\varphi (G_1)\subset G_2$.

It follows that if $(U,r_U,s_U)$ is a bi-submersion for  $\cF_1$ adapted to the path holonomy atlas of $\cF_1$ and $(V,r_V,s_V)$ is a bi-submersion for  $\cF_2$ adapted to the path holonomy atlas of $\cF_2$, then $U\circ V$ is a adapted to the path holonomy atlas of $\cF_2$.

\subsubsection{The morphism of foliation C*-algebras}

For every bi-submersion $(U,r,s)$ of $\cF_i$ ($i=1,2$), we denote by $\Omega^{1/2}U$ the bundle over $U$ of half densities on the bundle $\ker dr \oplus \ker ds$. 

Let $(U_1,r_1,s_1)$ be a bi-submersion for $\cF_1$ and $(U_2,r_2,s_2)$ a bi-submersion for $\cF_2$. Put $(W,r_W,s_W)=(U_1,r_1,s_1)\circ (U_2,r_2,s_2)$. Denote by $p_i:W\to U_i$ the projections. For $w=(u_1,u_2)\in W$, since $r_W=r_1\circ p_1$, and $\ker (dp_1)_w=\ker (dr_2)_{u_2})$, we find exact sequences $$0\to \ker (dr_2)_{u_2}\to \ker (dr_W)_{w}\to \ker (dr_1)_{u_1}\to 0$$ and in the same way, $$0\to \ker (ds_1)_{u_1}\to \ker (ds_W)_{w}\to \ker (ds_2)_{u_2}\to 0.$$ It follows that $\Omega^{1/2}W$ identifies with $\Omega^{1/2}U_1\otimes \Omega^{1/2}U_2$.

The algebra $C^*(M,\cF_i)$ is constructed in \cite[section 4.4]{AndrSk1}. Recall that for every bi-submersion $U$ in the path holonomy atlas of $\cF_i$ we have a linear map $\Theta^U_i:C_c^\infty (U;\Omega^{1/2})\to C^*(M,\cF_i)$, and that given an atlas $(U_j)_{j\in J}$ of bi-submersions adapted to the path holonomy atlas, the vector span of $\Theta^{U_j}_i:C_c^\infty (U_j;\Omega^{1/2}),\ j\in J$ is dense in $C^*(M,\cF_i)$.

\begin{proposition}\label{amidmopassant}
\begin{enumerate}
\item  There is a morphism $\Phi:C^*(M,\cF_1)\to \cM(C^*(M,\cF_2))$ characterized by the equality $\Phi(\Theta^{U_1}_1(f_1))\Theta^{U_2}_2(f_2)=\Theta^{U_1\circ U_2}_2(f_1\bullet f_2)$ for every bi-submersion $U_1$ for $\cF_1$, every bi-submersion $U_2$  for $\cF_2$, every $f_1\in C_c^\infty(U_1;\Omega^{1/2}U_1)$, $f_2\in C_c^\infty(U_2;\Omega^{1/2}U_2)$ where $f_1\bullet f_2\in C_c^{\infty}(U_1\circ U_2;\Omega^{1/2}(U_1\circ U_2))$ is the restriction to $U_1\circ U_2\subset U_1\times U_2$ of $f_1\otimes f_2\in C_c^{\infty}(U_1\times U_2;\Omega^{1/2}U_1\otimes  \Omega^{1/2}U_2)$.
\item The morphism $\Phi$ is non degenerate.
\end{enumerate}
\begin{proof}
\begin{enumerate}
\item We prove that, given a (faithful) representation $\varpi_2$ of $C^*(M,\cF_2)$, there is a representation $\varpi_1$ of $C^*(M,\cF_1)$ such that $\varpi_1(\Theta^{U_1}_1(f_1))\varpi_2(\Theta^{U_2}_1(f_2))=\varpi_2(\Theta^{U_1\circ U_2}_2(f_1\bullet f_2))$.

Recall from \cite[section 5]{AndrSk1} that the representation $\varpi_2$ corresponds to a triple $(\mu,H,\pi_2)$ where $\mu$ is a measure on $M$ which is quasi invariant by $G_2$, $H=(H_x)_{x\in M}$ is a $\mu$ measurable field of Hilbert spaces and $\pi_2$ is a representation of $G_2$ on the field $(\mu,H)$.

Fix a a bi-submersion $(U_1,r_1,s_1)$ for $\cF_1$.  Let $(U_2,r_2,s_2)$ be a bi-submersion for $\cF_2$, denote by $(W,r_W,s_W)$ the composition $U_1\circ U_2$ and let $p_i:W\to U_i$ be the projections. Put also $q=s_1\circ p_1=r_2\circ p_2$.

\begin{enumerate}
\item \emph{We show that $\mu $ is quasi-invariant by the groupoid $G_1$.} Choose positive smooth sections of the bundles of $1$-densities of the bundles  $\ker ds_1,\ \ker ds_2,\ \ker dr_1,\ \ker dr_2$. Thanks to these choices we construct measures $\mu\circ \lambda^{r_1}$ and $\mu\circ \lambda^{s_1}$ on $U_1$, measures $\mu\circ \lambda^{r_2}$ and $\mu\circ \lambda^{s_2}$ on $U_2$ and measures $\mu\circ \lambda^{r_W}$ and $\mu\circ \lambda^{s_W}$ on $W$. Consider also the measure $\mu \circ \lambda^q$. As $(W,r_W,s_W)$ is a bi-submersion for $\cF_2$ as well as $(W,q,s_W)$, it follows that the measures $\mu\circ \lambda^{r_W}$ and $\mu\circ \lambda^{q}$ are equivalent. Using a cover of $s_1(U_1)$ by open sets of the form $r(U_2)$, we deduce that $\mu\circ \lambda^{r_1}$ and $\mu\circ \lambda^{s_1}$ are equivalent. This being true for every bi-submersion $U_1$, it follows that the measure $\mu$ is quasi-invariant by the groupoid $G_1$.

\item \emph{We construct a representation $\pi_1$ of the groupoid $G_1$.} Associated with the bi-submersions $(U_2,r_2,s_2)$ and $(W,r_W,s_W)$ for $\cF_2$, we have measurable families of isomorphisms $\pi_{U_2}(u_2):H_{s_2(u_2)}\to H_{r_2(u_2)}$ and $\pi_{W}(u_1,u_2):H_{s_2(u_2)}\to H_{r_1(u_1)}$. We wish to define $\pi_{U_1}(u_1):H_{s_1(u_1)}\to H_{r_1(u_1)}$ by putting ({almost everywhere}) $\pi_{U_1}(u_1)=\pi_{W}(u_1,u_2)\circ \pi_{U_2}(u_2)^{-1}$. To do so, we define the measurable section $\tilde\pi_{U_1}$ over $W$ by setting $\tilde\pi_{U_1}(u_1,u_2)=\pi_{W}(u_1,u_2)\circ \pi_{U_2}(u_2)^{-1}$ and show it is constant ({almost everywhere}) along the fibers of $p_1:W\to U_1$. Let $(U'_2,r'_2,s'_2)$ be another bi-submersion for $\cF_2$. Since the identity of $U'_2\circ U_1\circ U_2$ is an isomorphism of  bi-submersions between  the composions $U'_2\circ (U_1\circ U_2)$ and $(U'_2\circ U_1)\circ U_2$, it follows that, ({almost everywhere} in $U'_2\circ U_1\circ U_2$) we have $$\pi_{U'_2}(u'_2)\pi_{U_1\circ U_2}(u_1,u_2)=\pi_{U'_2\circ U_1}(u'_2,u_1)\pi_{U_2}(u_2)$$  {whence}$$\pi_{U_1\circ U_2}(u_1,u_2)\circ \pi_{U_2}(u_2)^{-1}=\pi_{U'_2}(u'_2)^{-1}\circ \pi_{U'_2\circ U_1}(u'_2,u_1).$$

This shows that $\tilde\pi_{U_1}(u_1,u_2)$ is independent on $u_2$, \ie is constant along the fibers of $p_1$.
\end{enumerate}

The triple $(\mu,H,\pi_1)$ yields now a representation $\varpi_1$ of $C^*(M,\cF_1)$ on the Hilbert space $L^2(M,\mu ,H)$ of $L^2$ sections of the bundle $H$. 

Recall that the representations $\varpi_1$ and $\varpi_2$ are defined on any bi-submersion $(V,r_V,s_V)$ of $\cF_i$ by setting $\varpi_i\circ \Theta_i^V=\hat \pi_i^V$ where, for $f\in C_c^\infty(V;\Omega^{1/2}V)$ and $\xi\in L^2(M,\mu ,H)$, we have a formula: 
$$\hat\pi_i^V(f)(\xi)(x) = \int_{V^{x}}(1\otimes\rho^{V}(v)) (f(v))\pi_i^V(v)(\xi(s_{V}(v))) D^V(v)^{1/2}.$$
Recall that: \begin{itemize}
\item $V^x=r_V^{-1}(x)$;
\item $\rho^V$ is an isomorphism between $\Omega^{1/2}\ker ds_V$ and $\Omega^{1/2}\ker dr_V$ corresponding to the choice of a Riemannian metric on $M$, thus $(1\otimes\rho^{V}(v)) (f(v))$ is a $1$-density on $V^{x}$ and thus can be integrated;
\item the isomorphism $\rho^V$ is used to compare Lebesgue measures on the fibers of $s_V$ and of $r^V$; using it we obtain the Radon Nikodym derivative $D^V$ between the measures $\mu\circ \lambda^{s_V}$ and $\mu\circ \lambda^{r_V}$ which is used in this formula.
\end{itemize}
It is an easy check that, under the identification $\Omega^{1/2}(U_1\circ U_2)$ with $\Omega^{1/2}(U_1)\otimes  \Omega^{1/2}(U_2)$, we have $\rho^{U_1\circ U_2}(u_1,u_2)=\rho^{U_1}(u_1)\otimes \rho^{U_2}(u_2)$ and also $D^{U_1\circ U_2}(u_1,u_2)=D^{U_1}(u_1)D^{U_2}(u_2)$. 

One then checks immediately the equality $\hat \pi_1^{U_1}(f_1)\hat\pi_2^{U_2}(f_2)=\hat\pi^{U_1\circ U_2}_2(f_1\bullet f_2)$.

\item Take a bi-submersion $U_1$ for $\cF_1$, a bi-submersion $U_2$ for $\cF_2$. Let $(u_1,u_2)\in U_1\circ U_2$. If the class of $u_1$ is a unit of $G_1$ - \ie if $U_1$ carries an identity bissection through $u_1$, then the classes of $(u_1,u_2) $ and $u_2$ coincide. It follows that the bi-submersions of the form $U_1\circ U_2$ form  an atlas of bi-submersions equivalent to the path holonomy atlas of $\cF_2$ and therefore the span of $\Theta^{U_1\circ U_2}_2(f_1\bullet f_2)$ is dense in $C^*(M,\cF_2)$. The result follows. \qedhere
\end{enumerate}
\end{proof}
\end{proposition}

\subsection{Stability}

Since for every $x\in M$ there is $X\in \cF$ such that $X_x\ne 0$, we can find, using the construction of section \ref{Annexe}, \begin{itemize}
\item a locally finite open cover $(W_n)_{n\in \N}$;
\item for every $n\in \N$ a submanifold $U_n\subset W_n$ (closed in $W_n$) and a compact subset $U'_n\subset U_n$;
\item a vector field $X_n\in \cF$,
\end{itemize}
such that \begin{itemize}
\item $(t,u)\mapsto \exp_{tX_n}(u)$ is a diffeomorphism $f_n:{]}{-}1,1{[}\times U_n\to W_n$;
\item $\bigcup _{n\in \N}f_n({]}{-}1,1{[}\times U'_n)=M$.
\end{itemize}

It then follows from section \ref{lesVj}, that we can construct  a family of open sets $V_j$ satisfying the properties explained in step \ref{step2}.

The vector field $X_n$ defines a subfoliation $\cF_n=\{fX_n;\ f\in C_c^\infty(M)\}$ of $\cF$. Of course $X_n$ defines also an action $\beta_n$ of $\R$ on $M$ and $C^*(M;\cF_n)$ is a quotient of $C(M)\rtimes _{\beta_n}\R$ (the holonomy groupoid of $(M,\cF_n)$ is a quotient of $M\times _{\beta_n}\R$ since $M\times _{\beta_n}\R$ is an $s$-connected Lie groupoid with associated foliation $\cF_n$). Now it follows from the construction above that, if $n=n(j)$, the characteristic function $q_j$ of $V_j$ is a multiplier of $C(M)\rtimes _{\beta_n}\R$ and thus of $C^*(M;\cF_n)$ and finally, thanks to proposition \ref{amidmopassant}, of $C^*(M,\cF)$. The image of $q_j$ as a left multiplier of $C(M)\rtimes _{\beta_n}\R$ is $C_0(V_j)\Big(C(M)\rtimes _{\beta_n}\R\Big)$, therefore $q_jC^*(M,\cF)=C_0(V_j)C^*(M,\cF)$.

As $C^*(V_j,\cF_n)$ is stable and acts in a non degenerate way on $q_jC^*(M,\cF)$, it follows that the Hilbert-$C^*(M,\cF)$ module $q_jC^*(M,\cF)$ is stable. Finally, putting $q=\sum q_j$, we find that $qC^*(M,\cF)$ is stable.

To conclude, we prove that the Hilbert-$C^*(M,\cF)$ module $qC^*(M,\cF)$ is full. It follows from corollary \ref{corfinal} that the Hilbert-$C^*(M,\cF)$ module $C^*(M,\cF)$ is stable, \ie that the $C^*$-algebra $C^*(M,\cF)$ is stable.

\begin{proposition}
 Let $W\subset M$ be  open. Assume that $W$ meets all the leaves of $\cF$. Then $C_0(W)C^*(M,\cF)$ is a full Hilbert submodule of $C^*(M,\cF)$.
\begin{proof} 
Note that if $x,y\in M$ are in the same leaf, there is a bi-submersion $(U,r,s)$ in the path holonomy atlas and $u\in U$ such that $r(u)=x$ and $s(u)=y$. It follows that for every $x\in M$, there is a bi-submersion $(U,r,s)$ in the path holonomy atlas such that $x\in s(U)$ and $r(U)\subset W$. 

Let $\cU=(U_i,r_i,s_i)_{i\in I}$ be a family of bi-submersions defining the path holonomy atlas of $\cF$. For $i\in I$, put $U'_i=\{u\in U;\ r_i(u)\in W\}$. Define then $\widetilde \cU=(\tilde U_{i,j},\tilde r_{i,j},\tilde s_{i,j})_{(i,j)\in I^2}$ where $\tilde U_{i,j}=(U'_i)^{-1}\circ U'_j=\{(u,v)\in U_i\times U_j;\ r_i(u)=r_j(u)\in W\}$, and where we put $\tilde r_{i,j}(u,v)=s_i(u)$, and $\tilde s_{i,j}(u,v)=s_j(v)$.

Note that since $U'_j\circ \tilde U_{k,\ell}$ is adapted to $\cU$, it follows that $\tilde U_{i,j}\circ \tilde U_{k,\ell}$ is adapted to $\widetilde \cU$. We deduce that $\widetilde \cU$ is an atlas. It is obviously adapted to $\cU$; on the other hand, since $\cU$ is minimal, it is adapted to $\widetilde \cU$. 

It follows that the subspaces $C_c^\infty (\tilde U_{i,j})$ span a dense subalgebra of $C^*(M,\cF)$ and therefore $\{f^*g;\ f\in C^\infty_c(U'_i),\ g\in C^\infty_c(U'_j)\}$ span a dense subalgebra of $C^*(M,\cF)$. As the image of $C^\infty_c(U'_k)$ is in $C_0(W)C^*(M,F)$, the result follows.
\end{proof}
\end{proposition}

This ends the proof of theorem \ref{casdesfeuilletages}.

\begin{remarks}\begin{enumerate}
\item Note also that since the reduced $C^*$-algebra of the foliation is a quotient of the full one, it is also stable if $\cF$ vanishes nowhere.
\item As $C^*(M,\cF)$ corresponding to the minimal (path holonomy) atlas is a sub-algebra of the $C^*$-algebra of any atlas in a non degenerate way, it follows that, under our assumption that $\cF$ vanishes nowhere, all these $C^*$-algebras are stable.
\end{enumerate}
\end{remarks}

\selectlanguage{french}

\end{document}